\documentclass[a4paper,12pt,legno]{article}

\usepackage{amsmath}
\usepackage{amsfonts}
\usepackage{amsthm}
\usepackage{mathrsfs}
\usepackage{amssymb}

\newtheorem{thm}{THEOREM}
\newtheorem{lem}{LEMMA}
\newtheorem{prop}{PROPOSITION}
\newtheorem{cor}{COROLLARY}

\theoremstyle{remark}

\usepackage{graphicx}
\usepackage[T1]{fontenc}
\usepackage[polish]{babel}

\newcommand{\Rset}{\mathbb{R}}

\newcommand{\Zset}{\mathbb{Z}}

\begin{document}

\title{The strong mixing and the operator-selfdecomposability properties.}
\author{Richard C. Bradley\footnote{Department of Mathematics, Indiana University Bloomington, Indiana, USA.} \ and \ Zbigniew J. Jurek\footnote{Institute of Mathematics, University of Wroc\l aw, Wroc\l aw, Poland.} \,\footnote{Research funded by Narodowe Centrum Nauki (NCN) grant no Dec2011/01/B/ST1/01257.}}

\date{February 27, 2014}

\maketitle

\textbf{Abstract.}  For nonstationary, strongly mixing sequences of random variables taking their values in a 
finite-dimensional Euclidean space, with the partial sums 
being normalized via matrix multiplication, with certain
standard conditions being met, 
the possible limit distributions are precisely the operator-selfdecomposable laws.
\bigskip
\bigskip
\bigskip

   When one has observations (distributions) with values in an algebraic structure then their normalizations should be consistent with the structure in question. 
Thus, when $\xi_1,\xi_2, ...$ are $\Rset^d$-valued vectors then one should consider sums
\[
A_n(\xi_1+\xi_2+...+\xi_n)+x_n, \ \ \ \ \ \ \ \   (\ast)
\]
where $(A_n)$ are linear operators (matrices) on $\Rset^d$. Similarly,  if the algebraic structure is a Banach space or a topological group then the $(A_n)^{\prime}$s in $(\ast)$ should be bounded linear operators or automorphisms of the group, respectively. That novel paradigm required completely new algebraic methods and tools such as decomposability semigroups associated with probability measures,  
the Numakura Theorem on idempotents in (abstract) topological semigroups or elements of Lie theory. 
Sharpe (1969), for independent identically distributed 
$\xi_i^{\prime}$ s in $(\ast)$,  
and Urbanik (1972) and (1978), for infinitesimal triangular arrays $(A_n\xi_i, 1\le i \le n, n\ge 1)$, described limit distributions in the scheme $(\ast)$. 
The monograph Jurek and Mason (1993) summarized the research in that area for stochastically independent variables. However, the CLT for affine normalizations from 
Hahn and Klass (1981) still awaits for a coordinate-free proof.

   Here we will describe limiting distributions of $(\ast)$ for $\Rset^d$-valued random variables $\xi_1,\xi_2, ..., $ that are only strongly mixing, as defined by Rosenblatt (1956). Classical limiting distributions for strongly mixing sequences normalized by scalars are described in the monograph by Bradley (2007).

\medskip
\medskip
\textbf{1. Strong mixing and operator-selfdecomposability.}

\medskip
Let  $\Rset^d$ be  the  d-dimensional Euclidean space. As in
Jurek and Mason (1993), by $\bf{E}nd(\Rset^d) \equiv\bf{E}nd$ we denote the
Banach algebra of all bounded linear operators (matrices) on $\Rset
^d$ and by $\bf{A}ut(\Rset^d)\equiv \bf{A}ut$ the group of all
linear bounded and invertible operators (matrices). By
$\mathcal{P}(\Rset^d)\equiv \mathcal{P}$ we denote the topological semigroup of
all Borel probability measures on $\Rset^d$ with convolution $\ast$
and weak convergence topology.

Furthermore, let $(\Omega,\mathbb{F}, P)$ be a probability space reach
enough to carry uncountable family of independent uniformly
distributed random variables as well as sequences
$\textbf{X}:=(X_1,X_2,\dots)$ of $\Rset^d$- valued random vectors
(in short: random vectors); 
cf. Dudley (2002), Theorem 8.2.2.

We will say that a random vector $X$ or its probability distribution
$\mu$ is \emph{full or genuinely d dimensional } if its support is
not contained in any proper hyperplane in $\Rset^d$. (Recall
that\emph{a hyperplane} is  a linear subspace of $\Rset^d$ shifted
by a vector.) By $\mathcal{F}$ we denote the family of all full
measures. It is an open (in weak convergence topology) subsemigroup of $\mathcal{P}$ cf.
Jurek-Mason (1993), Corollary 2.1.2.

With a random vector $X$ or its  probability distribution $\mu$  we associate two semigroups of matrices:  the \emph{ Urbanik decomposability semigroup} $\mathbf{D}(X)$ (or  $\mathbf{D}(\mu)$) and the \emph{symmetry semigroup} $\mathbf{A}(X)$  (or
$\mathbf{A}(\mu)$) as follows: 
\begin{multline}
\mathbf{D}(X):=\{A\in \mathbf{End}:  X\stackrel{d}{=}AX+ Y \  \mbox{for some $Y$ independent of } \  X\}, \\
\mathbf{A}(X):=\{A\in \mathbf{End}:  X\stackrel{d}{=}AX+ a \  \mbox{for some vector } \  a \in \Rset^d \}, \ \ \ \ \ \ \ 
\end{multline}
where  $\stackrel{d}{=}$ denotes the equality in distribution. In an analogous way we define semigroups $\mathbf{D}(\mu)$ and $\mathbf{A}(\mu)$. Of course, $\mathbf{A}(X)\subset \mathbf{D}(X)$ and  the operators $0$ (zero) and $I$(identity) are always  in $\mathbf{D}(X)$.

   (The symbol $0$ will be used freely for the zero elements
of $\mathbb{R}$, $\mathbb{R}^d$, and $\mathbf{End}$.
In context, that should not cause confusion.) 

For the references below let us recall that 
\begin{multline}
(i)\ \ \  \mathbf{D}(\mu) \ 
\mbox{is a compact semigroup in}
 \  \  \mathbf{End}(\Rset^d) \ \  \  \ \mbox{iff} \  \ \  \  \mu \ \mbox{is a full measure} \\ 
 \ \ \  \mbox{iff}  \  \ \mathbf{A}(\mu)  \ \  \mbox{is a compact group in} \ \mathbf{Aut}(\Rset^d); \\
(ii)\ \ \  \mbox{If}\ \mu\ \mbox{is full, then}\ 
\mathbf{A}(\mu) \  \  \mbox{is the largest group in the Urbanik} \ \ \ \ \ \ \ \ \ \ \ \  \\
\mbox{semigroup}\ \mathbf{D}(\mu).\ \ \ \  
\end{multline}
Cf.\ Jurek and Mason (1993), Theorem 2.3.1 and 
Corollary 2.3.2 and Proposition 2.3.4.

We will say that \emph{a probability measure $\mu$ is operator-selfdecomposable} if there exist a sequence $b_n\in \Rset^d$,
a sequence $A_n\in \mathbf{End}$ and a sequence $X_n$ of 
\emph{independent} $\Rset^d$-valued random vectors such that

(i) the triangular array  $(A_nX_j: 1\le j \le n, n\ge 1)$  is infinitesimal;   

(ii) $\lim_{n\to \infty}  A_n(X_1+X_2+...+X_n)+b_n 
\Rightarrow \mu$.

\noindent (Condition (i) simply means that
$A_nX_j \to 0$ in probability as $n \to \infty,\
j \in \{1,\dots, n\}$.)

The main characterization  due to K. Urbanik is as follows: 

\emph{A full measure $\mu$ is operator-selfdecomposable iff its decomposability semigroup}   $\textbf{D}(\mu)$ 
\emph{contains at least one one-parameter semigroup} 
$\{\exp (-tQ), t \geq 0\}$ \emph{with} 
$\exp (-tQ) \to 0$ (the zero matrix) \emph{as} $t \to \infty$. 

Cf.\ Jurek-Mason (1993), Theorem 3.3.5.
(The stipulation in that theorem that $Q$ be invertible,
is superfluous; note that $Q^{-1}=\int_0^{\infty}e^{-sQ}ds$; the integral  is well defined because  $\exp (-tQ) \to 0$ as $t \to \infty$.) 

Also it might be of some importance to mention here
that we have  the following random integral representation:
\[
\mbox{ $\mu$ is operator-selfdecomposable iff} \ \  \mu=\mathcal{L}\big(\int_0^{\infty}e^{-tQ}dY(t),
\big)
\]
 for some L\'evy process $Y$ ( so called \emph{background driving L\'evy process } (BDLP); cf.  Jurek (1982) or Jurek-Mason (1993), Theorem 3.6.6.

\medskip
Since our aim here is to extend the notion of operator-selfdecomposablity to some dependent random variables, 
let's recall that for two sub-$\sigma$-fields 
$\mathcal{A}$ and $\mathcal{B}$  of
$\mathbb{F}$ we define \emph{the measure of dependence} 
$\alpha$ between
them as follows:
\[
\alpha(\mathcal{A},\,\mathcal{B}):=\sup_{A\in\mathcal{A}, B\in
\mathcal{B}}\,|P(A\cap B)-P(A)P(B)|\,.
\]
For a given sequence $\textbf{X}:=(X_1,X_2,\dots)$ of 
$\mathbb{R}^d$-valued
random variables, we define for each positive integer $n$ the
dependence coefficient
\begin{equation}
\alpha(n)\equiv \alpha(\textbf{X};n):=\sup_{j\in {\bf
N}}\,\alpha\big(\,\sigma(X_k, 1\le k \le j), \, \, \sigma(X_k, k \ge
j+n)\big),
\end{equation}
where $\sigma(\dots)$ denotes the $\sigma$-field generated by
(\dots). We will say that a sequence \textbf{X} is \emph{strongly
mixing} (Rosenblatt (1956)) if
\begin{equation}
\alpha(n)\to 0 \ \ \mbox{as} \ \ n\to \infty.
\end{equation}
Of course, if the elements of \textbf{X} are stochastically
independent then $\alpha(\textbf{X};n)\equiv 0$.

\begin{thm}
 Let $\textbf{X}:=(X_1,X_2,...)$ be sequence $\Rset^d$-valued random vectors with the partial sums $S_n:=X_1+X_2+...+X_n$,
 and let $(A_n)\in \bf{E}nd(\Rset^d)$ and $(b_n)\in\Rset^d$ be sequences of bounded linear operators and vectors, respectively
 satisfying conditions:

(i) $\alpha(\textbf{X};n)\to 0$ as $n\to \infty$,  i.e., the
sequence $\textbf{X}$ is strongly mixing;

(ii) the triangular array $(A_n\,X_j, \ 1 \le j \le n, n\ge 1)$ is
infinitesimal;

(iii) $A_nS_n+b_n\Rightarrow \mu$, for some full probability measure
$\mu$.

\noindent Then the limit distribution $\mu$ is
operator-selfdecomposable, that is, there exists a one parameter semigroup 
$\{e^{-tQ}: t\ge 0\}\subset \mathbf{D}(\mu)$ 
with 
$\lim_{t\to \infty}e^{-tQ}=0$.
\end{thm}
\bigskip

The line of reasoning in our proof of this theorem
is as follows. First, in Section 2 we investigate the normalizing sequence $(A_n)$ of matrices, and in 
particular we show that one may choose a more appropriate sequence $(\widetilde{A}_n)$. 
Then in Section 3, using the new normalizing sequence, 
we construct in a few steps 
a one-parameter semigroup $(e^{tQ}, t\ge 0)$. 
Here, we follow Urbanik (1972); but we could also argue similarly as in Urbanik (1978) or Jurek-Mason (1993), where the proof is valid in infinite dimensional linear spaces.

\medskip
\medskip

\vspace{0.25in} \textbf{2. Auxiliary propositions and lemmas.}

\medskip
First,  some consequences  of  the operator-convergence of types
theorems   (Section 2.2 in Jurek -Mason (1993)):

\begin{prop}
Under the assumptions (ii) and (iii) in Theorem 1,

\noindent a) $A_n \to 0$ as $n \to \infty$; the inverse
$A_n^{-1}$ {exists for all sufficiently large $n$;

\noindent
b) there exist $\widetilde{A}_n$ for which (ii) and (iii) 
(in Theorem 1) hold and
 
$ \widetilde{A}_{n+1}\widetilde{A}_n^{-1}\to I $
(the identity matrix)};

\noindent c) one has that
\[  
\lim_{n\to \infty}\,
\Bigl|\frac{det A_{n+1}}{det A_n}\Bigl|
=\lim_{n\to
\infty} \frac{det \widetilde{A}_{n+1}}
{det \widetilde{A}_n} =1.
\qquad \qquad \qquad
\]
\end{prop}

Cf. Jurek-Mason (1993), Section 3.2 : Propositions 3.2.1 and  3.2.2. For part c) one needs also Corollaries  2.3.2 and 2.4.2 as
$\widetilde{A}_n:= H_nA_n$  for some $H_n\in\mathbf{A}(\mu)$.

\medskip
Second, a note on uniform infinitesimal triangular arrays. 

\begin{lem}
Suppose that for each $n\in \mathbb{N}$, $I_n$ is a nonempty set.
Suppose that for each $n\in\mathbb{N}$ and each $j\in\ I_n$,
$X_{n,j}$ is Banach space valued random element. The the following
two statements are equivalent:
\begin{equation}
(A) \qquad \qquad \forall(\epsilon>0) \  \ \lim_{n\to
\infty}\sup_{j\in I_n}\,P(||X_{n,j}||\ge\epsilon)=0
\end{equation}
There exists a sequence $\delta_1 \ge \delta_2 \ge ... \ge
\delta_{n-1} \ge \delta_{n} \to 0$ as $n\to\infty$, such that
\begin{equation}
(B) \qquad \quad \forall (n\in \mathbb{N})\quad \forall (j\in I_n)
\quad P(||X_{n,j}||\ge \delta_n)\le \delta_n
\end{equation}
\end{lem}
\emph{Proof.}  ($(A) \Rightarrow (B))$ Let $N_1:=1$. For already
defined $N_1, N_2,...,N_{m-1}$, let $N_m>N_{m-1}$ be such that
\begin{equation}
\forall ( n\ge N_m) \quad \forall(j\in I_n) \quad P(||X_{n,j}||\ge
1/m)\le \frac{1}{m}
\end{equation}
which always exists by (5). Next, for the defined sequence
\[
1=N_1<N_2<...<N_n<...,
\]
let us define the sequence $(\delta_n)$ as follows:
For each $m \in \mathbb{N}$, 
\begin{equation}
\delta_n:=1/m, \quad \mbox{for all $n$ such that} \quad   N_m\le
n\le N_{m+1}-1 .
\end{equation}
Thus by virtue of the above construction,
$\delta_n \to 0$ as $n \to \infty$; 
and for each 
$n\in\mathbb{N}$ there exists exactly one
$m\in\mathbb{N}$ such that $N_m\le n\le N_{m+1}-1$,
and by (7) and (8),
\[
P(||X_{n,j}||\ge \delta_n)=P(||X_{n,j}||\ge 1/m)\le 1/m=\delta_n  \ \  \mbox{for all} \  j \in I_n,
\]
which completes the proof $(A)\Rightarrow (B)$. The implication 
$(B)\Rightarrow (A)$ is obvious.

\begin{cor}
For the  infinitesimal triangular  array $(X_{n,j})$ as in Lemma 1 and $q_n\to \infty$, $q_n\le \delta_n^{-1/2}$ then for any set $Q\subset  I_n$  such that  $card  Q\le q_n$
we have that
\[
P(||\sum_{k \in Q} X_{n,k}||\ge \delta_n^{1/2})\le \delta_n^{1/2} 
\]
\end{cor}

Since $ card\, Q\le q_n$  and $\delta_n\,q_n\le \delta_n^{1/2}$
\begin{multline*}
P(|| \sum_{k\in Q}X_{n,k} ||\ge \delta_n^{1/2})  \le P(||\sum_{k\in Q} X_{n,k}||\ge \delta_nq_n) \\   \le \sum_{k\in Q}\, P(||X_{n,k} ||\ge \delta_n)   
\le q_n\delta_n \le \delta_n^{1/2}.
\end{multline*}

\medskip
Third,  a generalization of Proposition 3.2.3 in Jurek-Mason, for strongly mixing sequences.
\begin{prop}
Suppose the hypothesis of Theorem 1, including all of
conditions (i), (ii), and (iii) there, hold.
Suppose also that for every $n \in \mathbb{N}$, the
matrix $A_n$ is invertible. 
Then

\begin{equation} 
\sup\{||A_nA_m^{-1}||:  
 n \in \mathbb{N}, 1\le m \le n\}<\infty. 
\end{equation}

Moreover, if for each $n \in \mathbb{N}$,
$m_n$ is an integer such that $1 \le m_n \le n$, then 
all limits points of the sequence 
$(A_nA_{m_n}^{-1})_{n\in \mathbb{N}}$  
are in $ \textbf{D}(\mu)$.
\end{prop}

\emph{Proof.} We shall first prove (9).  
Suppose that for each $n \in \mathbb{N}$, $m_n$ is an
integer such that $1 \le m_n \le n$.
To prove (9), it suffices to prove that 
\begin{equation}
\sup_{n \in \mathbb{N}}||A_nA_{m_n}^{-1}|| < \infty.
\end{equation}
If instead $||A_nA_{m_n}^{-1}|| \to \infty$ along
some subsequence of $n \in \mathbb{N}$, then within
that subsequence the integers $m_n$ could not be
bounded (for otherwise $||A_nA_{m_n}^{-1}|| \to 0$
would occur along that subsequence by 
Proposition 1(a)), and there would be a further
subsequence along which $m_n \to \infty$.
Letting $m_n := n$ for all $n$ not in that
``further subsequence,'' we have reduced our task
(for the proof of (9)) to proving (10) under the
additional assumption that $m_n \to \infty$
as $n \to \infty$.    

\medskip
For the rest of the proof it is assumed that
\[ m_n\to \infty,  \ \   0<q_n\to \infty \ \ \mbox{and} \ \ \  q_n \le \delta_n^{ -1/2},
\]
where the sequence  $(\delta_n)_{n\in \mathbb{N}}$  
is as in Lemma 1.

\medskip
For $n \in \mathbb{N}$ define random vectors as follows:
\[
\eta_n: =
\begin{cases}
0,  & \ \ \ \mbox{if} \ \ m_n=n; \\
S(\textbf{X}, n)-S(\textbf{X}, m_n),  & \ \ \ \mbox{if} \ \ \ n- q_n\le m_n<n \\
S(\textbf{X}, m_n+q_n)-S(\textbf{X}, m_n),  & \  \  \ \mbox{if} \ \ \ m_n\le  n-q_n-1.
\end{cases}
\]
and
\[
\xi_n:= 
\begin{cases}
0,& \ \ \ \mbox{if} \ \ \ n-q_n \le m_n\le n  \\
S(\textbf{X}, n)-S(\textbf{X},  m_n +q_n), & \ \ \ \mbox{if} \ \ \  m_n \le n-q_n -1.
\end{cases}
\]
Thus
\begin{equation}
S_n\equiv S(\textbf{X},n)= S(\textbf{X}, m_n)+\eta_n+\xi_n
\end{equation}
Since $\eta_n$ is either zero or the sum of 
at most $q_n\le \delta_n^{-1/2}$ of the variables $X_k$,
one has by Corollary 1 that 
$A_n\eta_n\to 0 $ in probability as $n\to \infty$. Consequently,
\begin{equation}
A_nS(\textbf{X}, m_n)+ A_n\xi_n +b_n \Rightarrow Z , \ \ \mbox{as} \ \ \ n \to \infty.
\end{equation}
>From the description of $\xi_n$'s, for the case $m_n\le n-q_n-1$ we have
\begin{equation}
\alpha(\sigma(S(\textbf{X}, m_n)),\sigma (\xi_n))\le \alpha (\textbf{X}, q_n+1) 
\end{equation}
In the opposite case ($m_n\ge n-q_n$), the $\xi_n$'s are 
zero (constant variables) and the left-hand side of (13)  
is therefore zero.
Now from (12),
\begin{equation}
\big[A_nA_{m_n}^{-1}\big(A_{m_n}S(\textbf{X}, m_n)+b_{m_n}\big)\big] + \big[A_n\xi_n+b_n-A_nA^{-1}_{m_n}b_{m_n}\big] \Rightarrow Z, n\to\infty
\end{equation}
For simplicity, let $V_n$ and $W_n$ denote the first and the second expressions in the above square brackets, that is
\begin{equation}
V_n+W_n \Rightarrow Z
\end{equation}
>From (13) and Corollary 1.11 in Bradley (2007), Vol.\ 1,
\begin{multline*}
|\mathbb{E}[\exp\, i<t,V_n+W_n>]-\mathbb{E}[\exp\, i<t,V_n>]\cdot \mathbb{E}[\exp\,i<t,W_n>]| \\
\le 16\, \alpha(\mathbf{X}, q_n+1)\to 0  \ \ \ \mbox{as} \ \ \ n\to \infty.
\end{multline*}
Hence by (15),
\begin{equation*} 
\mathbb{E}[\exp\, i<t,V_n>]\cdot \mathbb{E}[\exp\,i<t,W_n>]\to \mathbb{E}[\exp i<t, Z>]\  \mbox{as} \  n\to \infty.
\end{equation*}
Our next task is to replace vectors  $W_n$  by vectors that are stochastically independent of $V_n$. To this aim, let
$\zeta_{n, m}\, ,  m=1,2,..,n$ be random vectors independent of $\sigma (\mathbf{X}, Z)$ such that
\begin{multline*}
\zeta_{n,m}:=b_n-A_nA^{-1}_m\,b_m, \ \ \mbox{if} \ \ \ m+q_n\ge n;  \  (\mbox{a constant}); \ \mbox{for} \ m+q_n<n,  \\
\mathcal{L}(\zeta_{n,m}):= \mathcal{L}(A_n\big(S(\textbf{X},n)- S(\textbf{X}, m+q_n))+b_n-A_nA_m^{-1} b_m\big). 
\end{multline*}
But note that   $\mathcal{L}(W_n)=\mathcal{L}(\zeta_{n,m_{n}})$ and thus by the above,
\[ \mathbb{E}[\exp\, i<t,V_n>]\cdot \mathbb{E}[\exp\,i<t, \zeta_{n,m_n}>]\to \mathbb{E}[\exp i<t, Z>]\  \mbox{as} \  n\to \infty
\]
and therefore
\[
\mathcal{L}(V_n)\ast\mathcal{L}(\zeta_{n,m_n})=\mathcal{L}(V_n+\zeta_{n,m_n})\Rightarrow  Z
\]
By Parthasarathy (1967), Theorem 2.2 in Chapter III , $(\mathcal{L}(V_n))_n$ is shift compact, or 
in the ``symmetrization'' terminology there,
\[
(A_nA_{m_n}^{-1}\big(A_{m_n}S(\textbf{X}^{\circ}, m_n)\big))_{n \in \mathbb{N}}, \ \  \mbox{is compact and} \ \  A_{m_n}S(\textbf{X}^{\circ}, m_n)\Rightarrow Z^{\circ}.
\]
>From Lemma 2.2.3 in Jurek-Mason (1993) we now get the boundedness 
in Proposition 2. 

Further, if $D$ is a limit point of the family of matrices  
$(A_nA_{m_n}^{-1})_{n \in \mathbb{N}}$ then from (14) 
we get 
$\mathcal{L}(DZ+Y)=\mathcal{L}(Z)$ for some random
variable $Y$, a limit 
point of $(\zeta_{n,m_n})_{n\in\mathbb{N}}$, 
independent of Z. 
This completes the proof Proposition 2.

\medskip
\medskip
\medskip
\textbf{ 3. Construction of the one-parameter semigroup in $\mathbf{D}(\mu)$.}

\medskip
Here we follow the Urbanik construction from Urbanik (1972); see also Jurek-Mason (1993), Section 3.3.
Throughout this section, as in the hypothesis of Theorem 1,
we assume that {\it the probability measure $\mu$ is full\/};
and as allowed by Proposition 1, we assume that the 
matrices $A_n$ satisfying conditions (ii) and (iii) 
in Theorem 1 are invertible and satisfy 
$A_{n+1} A_n^{-1} \to I$ 
and (hence) $det A_{n+1}/det A_n \to 1$ 
(as $n \to \infty$). 

\medskip
By \emph{an idempotent $J$}, in $\mathbf{End}(\Rset^d)$, we mean a projector from $\Rset^d$ onto  the linear subspace $J(\Rset^d)$, that is $J^2=J$.  
Following Numakura (1952), p. 103, we will say that \emph{idempotent K is under idempotent J}, if  
$K \neq J$ and $JK=KJ=K$.  Hence, in particular,  $K(R^d)\subsetneq J(R^d)$.
If there is no non-zero idempotent under J, the we will say that J is \emph{a primitive idempotent}.   

Idempotents will play a crucial role below as we have the following: \emph{ for an idempotent $J$ we have that}
\[
J\in\mathbf{D}(\mu) \ \ \ \mbox{iff} \ \ \  
(I-J)\in\mathbf{D}(\mu) \ \  \ \mbox{and} \ \  \mu=J\mu\ast(I-J)\mu
\]
Furthermore, if an idempotent $K$ is under $J$ and both are in $\mathbf{D}(\mu)$ then
\begin{equation}
\mu = K\mu\ast(J-K)\mu\ast(I-J) \mu \ \  \mbox{and} \ \ K+(J-K)+(I-J)=I
\end{equation}
for details cf. Jurek-Mason (1993), Theorem 2.3.6.

Below $\det_JA$ means the determinant of matrix representation of the operator $JA$ in $J(\Rset^d)$ relatively to an orthogonal basis of $J(\Rset^d)$. Hence
we get
\begin{multline}
(a) \ \ det_JA=det_J(JA)= det_J (AJ)=det_J(JAJ); \\ 
(b) \ \ det_J(AJB)=det_JA\,det_JB ;   \quad  \quad \quad \quad
\qquad \qquad \qquad \qquad \\  
(c)\ \  det (JAJ+(I-J)B(I-J))=det_J A\,det_{I-J} B.
\qquad 
\end{multline}

\medskip
\begin{lem} 
For a given  idempotent $J\in \mathbf{D}(\mu)$, for each $0<c<1$ there exist $K_c\in\mathbf{D}(\mu)$ such that
$\det_JK_c=c$
\end{lem}
\emph{Proof.}  For $1 \leq n \leq m$, one has the 
inequalities
\begin{multline*}
 ||A_{m+1}A_n^{-1}-A_mA_n^{-1}|| \le || A_mA_n^{-1}||\,  ||A_{m+1}A_m^{-1} - I|| \\
 \le ( \sup_{n\le m}||A_mA_n^{-1}||)\,||A_{m+1}A_m^{-1}-I||.
\end{multline*}
 Since by Proposition 2, $\sup\{||A_nA_m^{-1}||:1\le m\le n, n\in \mathbb{N} \}<\infty, $  we get
\[
\lim_{m\to \infty} \sup_{n\le m} |\, ||A_{m+1}A_n^{-1}||-||A_mA_n^{-1}||\,|\le \lim_{m\to \infty} \sup_{n\le m}  ||A_{m+1}A_n^{-1}-A_mA_n^{-1}||=0.
\]
Since the functions $\Rset^{d^2}\ni A\to ||A||$ and $ \Rset^{d^2}\ni A\to \det_J A $ are continuous therefore by putting $b_{m,n}:=\det_JA_mA_n^{-1}$  $(n\le m)$) we infer that
\begin{equation}
b_{n,n}=1,  \ \ \ \lim_{m\to \infty} b_{m,n}=0  \ (n=1,2,...), \ \ \lim_{m\to \infty} \sup_{n\le m}|b_{m+1,n}-b_{m,n}|=0.
\end{equation} 
Thus  for any $0<c<1$  and the sequence $m_n:=\sup\{k\ge n: b_{k,n}\ge c\}$ we get   $b_{m_n+1,n } <c \le  b_{m_n,n}$,  so  from (18), $\lim_{n\to \infty}b_{m_n,n}=c$.  

Furthermore, by Proposition 2,  if $K_c$ is a limit point of a sequence $(A_{m_n} A_n^{-1})$ then $K_c$  is in $\mathbf{D}(\mu)$  and, by (18), $\det_J K_c=c$, which concludes the proof.

\begin{lem}
Let $J$ be non-zero idempotent in $\mathbf{D}(\mu)$. Then there exists $T_n\in\mathbf{D}(\mu), n=1,2,...$ such that
\begin{equation}
 JT_n=T_nJ=T_n, \ \ \    T_n\to J \  \mbox{and} \  \lim_{k\to \infty}T_n^k=0 \ \  (n=1,2,...)
\end{equation}
\end{lem}
\emph{Proof.} We shall justify the above claim by the mathematical  induction with respect to the dimension of linear space $J(\Rset^d)$. 

\medskip
Step 1.  $\dim J(\Rset^d)=1$.

>From Lemma 2, there exist 
$K_n \in \mathbf{D}(\mu)$ such that $det_J K_n=1-1/n$. 
Putting $T_n:=JK_nJ$ we have that  the linear  transformation $T_n: J(\Rset^d)\to J(\Rset^d)$  must be  a multiple of $J$; 
($ \dim J(\Rset^d)=1)$. But $\det_J T_n=\det_J K_n=1-1/n$ and thus $T_n=(1-1/n)J$ which satisfies (19).

\medskip
Step 2. Assume $\dim J(\Rset^d)= l >1$ and for all idempotents $ K\in\mathbf{D}(\mu)$ such that $\dim K(\Rset^d)<l$, Lemma 3 is true.

\medskip
Case (i).  Assume that there exist non-zero idempotent $L\in \mathbf{D}(\mu)$ such that $L\neq J$  and  
\begin{equation}
L=JL=LJ, 
\end{equation}
that is, the idemptent J is not a primitive one.

 From the above $J-L$ is also an idempotent. From Jurek-Mason (1993) Theorem 2.3.6 (a) ,  $I-L\in \mathbf{D}(\mu)$. Hence $J(I-L)=J-L\in \mathbf{D}(\mu)$ .

Since $\dim L(\Rset^d)<l$ and $\dim (J-L)(\Rset^d)<l$ therefore, by the mathematical induction assumption, there exist  sequences $(U_n)$ and $(V_n)$  in $\mathbf{D}(\mu)$ such that
\begin{multline*}
 U_n\to L,  \ \ LU_n=U_nL=U_n   \ \  \mbox{and} \  \lim_{k\to \infty}U_n^k=0 \ \  (n=1,2,...) \\
 V_n\to J-L,  \ \ (J-L)V_n=V_n(J-L)=V_n   \ \  \mbox{and} \  \lim_{k\to \infty}V_n^k=0 \ \  (n=1,2,...)
\end{multline*}
Then putting $T_n:=U_n+V_n$ we have $T_n \to J$. Further, 
from the identity $T_n=LU_nL+(I-L)V_n(I-L)$ and  again by Theorem 2.3.6 (d) in Jurek-Mason (1993)
we get that $T_n\in\mathbf{D}(\mu)$ and also  $T_n^k=U_n^k+V_n^k\to 0$ as $k\to \infty$. 
(Also, the first two equalities in (19) hold by
an elementary argument.)\ \  
This completes the Case (i).

\medskip
Case(ii). There are no non-zero idempotents $L$  in $\mathbf{D}(\mu)$ different from $J$ and satisfying  $JL=LJ=L$, i.e.,
idempotent J is a primitive idempotent.

\medskip
>From Lemma 2, choose $D_n\in \mathbf{D}(\mu)$ such that
\begin{equation}
0< det_J D_n<1  \ \ \ \mbox{and}  \ \ \  \lim_{n\to \infty} det_J D_n=1.
\end{equation}
By (a) in the formula (17) 
we may assume that $JD_n=D_nJ=D_n$ and if $D$ is a limit point of the sequence $D_n$ then  we also have equalities
\begin{equation}
D\in\mathbf{D}(\mu), \ \ \  JD=DJ=D  \ \ \  \  \mbox{and}  \ \  det_{J} D=1
\end{equation}
Put $A:=D+I-J$. Note that $A=JDJ+(I-J)I(I-J)$. Then  by (d) in Theorem 2.3.6 from Jurek-Mason (1993) $A\in\mathbf{D}(\mu)$. However,  
by (17) and (21), 
\[
 det A=det_{J+I-J}(JDJ+(I-J)I(I-J))=det_J D\,det_{I-J}(I-J)=1.
\]
Consequently, by Jurek-Mason (1993), Proposition 2.3.5 and
Corollary 2.3.2, 
\begin{equation}
A\in \mathbf{A}(\mu) \  \mbox{(a compact group in} \textbf{A}ut) \  \mbox{and} \   A^{r_n}\to I,
\end{equation}
for some $r_1<r_2<\dots$. Since $JA^n=D^n$ we have that $D^{r_n}\to J$.  Furthermore, since $ D$ is a limit point of the sequence $D_n$ we can choose a subsequence $(k_n)$ such that
\[
T_n:=D^{r_n}_{k_n}\to J ;  \ \ JT_n=T_nJ=T_n; \ \ 0<det_J T_n<1.
\] 
To complete the proof one needs to show that $T_n^k\to 0$ as $k\to \infty$.

For each $n$, the monothetic  semigroup $sem (T_n)$  (the smallest closed subsemigroup containing  $T_n$)  is compact in $\mathbf{D}(\mu)$. By the Numakura Theorem (Corollary 1.1.3  in Jurek-Mason) the limit points of $(T_n^k)_{k\in \mathbb{N}}$ form a group, denoted  by $K(T_n)$, with the unit $L$ that satisfies
\[
JL=LJ=L \ \  \mbox{and}  \ \  det_JL=0 \ \ \mbox{and thus} \ \ L\neq J 
\]
Because of the assumption (ii) we must have $L=0$. Consequently $T_n^k\to 0$ as $k\to \infty$, which completes the proof of Lemma 3.

\medskip
Using the formula (16) inductively, there are finitely many non-zero  primitive idempotents $J_1,J_2,..., J_q$ in $\mathbf{D}(\mu)$ , $q\le d$ (the dimension of $\Rset^d$),  such that
\begin{equation}
 I=J_1+J_2+...+J_q, \ \ \  J_rJ_s=J_sJ_r=0 \  (1\le r \neq s \le q). 
\end{equation}
Thus, in particular, for every $s$ there is no non-zero idempotent $K$ such that
\begin{equation} 
J_sK=KJ_s=K.
\end{equation}
Finally, recall that for idempotents satisfying (24) (not necessary primitive ones) we have
\begin{equation}
\mbox{\emph{if}} \ A_1,A_2,...,A_q\in\mathbf{D}(\mu)  \ \ \mbox{\emph{then}} \ \ J_1A_1J_1+J_2A_2J_2+ ...+J_qA_qJ_q\in \mathbf{D}(\mu);
\end{equation}
for details cf. Jurek-Mason (1993), Theorem 2.3.6.

\begin{lem}
There exists a positive integer $q$ and 
a one parameter semigroup $\{C_w:  w \in W \}\subset \mathbf{D}(\mu)$ (where $W$ denotes the set of
non-negative rational numbers) such that $\det C_w=e^{-qw}$
and $C_0 = I$.
\end{lem}
\emph{Proof.}
In view of Lemma 3, for the idempotents $J_r$ in (24),  let us choose $T_{n,r}\in \mathbf{D}(\mu)$  such that for $1\le r\le q, n\ge 1$ we have
\begin{equation}
 J_rT_{n,r}=T_{n,r}J_r=T_{n,r}, \ \ \    T_{n,r} \to J_r, \  \ \  \lim_{k\to \infty}T_{n,r}^k=0, \ \ \  0< det_{J_{r}} T_{n,r}<1.
\end{equation}
Note that $\lim_{n\to\infty} (\log det_{J_r} T_{n,r})= 0$  and put  
$d(n,r):= [ (- \log det_{J_r} T_{n,r})^{-1}]$, where the bracket $[.]$ denotes the integer part. Hence,
$\lim_{n \to \infty} d(n,r) = \infty$ and
\begin{equation}
\lim_{n\to \infty}(d(n,r)\cdot(-\log det_{J_r}T_{n,r}))=1,  \quad (r=1,2,...,q). 
\end{equation}
Further, let $W$ denote the set of all non-negative rational numbers (as in the statement of Lemma 4).  Then
\begin{multline}
 T_{n,r}^{[w\,d(n,r)]}\in \mathbf{D}(\mu) \ \ \mbox{for all} \  n\in\mathbb{N}, w\in W,  1\le r\le q;  \  \ \ \mbox{and by (26)}, \\ \ \  \ \ \
\sum_{r=1}^q J_{r} T_{n,r}^{[w\,d(n,r)]}J_{r} \in\mathbf{D}(\mu)  \ \ \quad  \mbox{for all} \  w\in W, \quad  n\in\mathbb{N};  
\end{multline}
Since $\mathbf{D}(\mu)$ is compact, there exist a subsequence $Q\subset \mathbb{N}$ and $C_w\in\mathbf{D}(\mu)$ such that for each $w\in W$
\begin{equation}
\sum_{r=1}^q  T_{n,r}^{[w\,d(n,r)]}= \sum_{r=1}^q J_{r} T_{n,r}^{[w\,d(n,r)]}J_{r}\to C_w, \ \mbox{as} \ \  n\to \infty, n\in Q.
\end{equation}
(Note for $w=0$ that this gives $C_0 = I$ by (24).)\ \ 
Hence, from (28)  and (30) we get
\begin{equation}
 det_{J_s} C_w=\lim_{n\to \infty, n\in Q} det_{J_s}T_{n,s}^{[w\,d(n,s)]}=\lim_{n\to \infty, n\in Q} (det_{J_s}T_{n,s})^{[w\,d(n,s)]}= e^{-w}.
\end{equation}
So, by (17) and (24) we conclude
\begin{multline}
det\,C_w= \lim_{n\to \infty}det_{J_1+...+J_q} \big(\sum_{r=1}^q J_{r} T_{n,r}^{[w\,d(n,r)]}J_{r}\big)  \\ = \prod_{r=1}^q\lim_{n\to \infty, n\in Q} det_{J_r}T_{n,r}^{[w\,d(n,r)]}=  e^{-qw}
\end{multline}
 To show that $\{C_w: w\in W\}$ is indeed  a one-parameter additive semigroup, note that for the integer 
part function $a\ni \Rset\to [a]\in\Zset$ (integers)  we have
\begin{equation} 
 [a+b]-[a]-[b]\in \{0,1\}
\end{equation}
(because  $a+b-1<[a+b]\le a+b$, \ \  $ -a\le -[a]<1-a$\  and \ $-b\le -[b]<1-b$)

Hence, for $w\in W$ and $u\in W$,
\[
s_n:=[(w+u)\,d(n,r)]-[w\,d(n,r)]-[u\,d(n,r)]\in \{0,1\},    \ (r=1,2,..,q)
\]
Hence by (27)
\begin{equation}
\lim_{n\to \infty}J_{r}\big( T_{n,r}^{s_n}- I\big)J_{r} =0, \  (r=1,2,...,q),
\end{equation}
since $s_n=0$ or $s_n=1$.
Finally, from (30) and (24),
\begin{multline*}
C_{w+u}-C_wC_u = \\   \lim_{n\to\infty}  \sum_{r=1}^q J_{r} T_{n,r}^{[(w+u)\,d(n,r)]}J_r -\big(\lim_{n\to \infty}  \sum_{r=1}^q J_{r} T_{n,r}^{[w\,d(n,r)]}J_r\big)\big( \lim_{n\to \infty}  \sum_{s=1}^q J_{s} T_{n,s}^{[u\,d(n,s)]}J_s\big)
\\ = \sum_{r=1}^q\lim_{n\to\infty} T_{n,r}^{[w\,d(n,r)]+[u\,d(n,r)]}J_{r}\big( T_{n,r}^{s_n}- I\big)J_{r}.  \ \quad \qquad
\end{multline*}
Since $T_{n,r}^k\in\mathbf{D}(\mu) \ (n,k\in\mathbb{N}, r=1,2,...,q)$  and $\mathbf{D}(\mu)$ is compact (thus 
the norms of its members are bounded, say by B), one has 
from above and (34),
\[
||C_{w+u}-C_wC_u||\le B\sum_{r=1}^q\lim_{n\to \infty}||J_{r}\big( T_{n,r}^{s_n}- I\big)J_{r}  ||=0,
\]
which gives the one-parameter semigroup property  $C_{w+u}=C_wC_u$.

\medskip
\begin{lem}
For the given (full) probability measure $\mu$, its Urbanik decomposability semigroup $\mathbf{D}(\mu)$ contains at least one one-parameter semigroup $\{e^{-tQ}, t \geq 0\}$ ($Q$ is a matrix)  such that $e^{-tQ}\to 0$, as $t\to \infty$.
\end{lem}
\emph{Proof.}  Throughout this proof, we use freely
all notations and arguments in the {\it proof\/} (as well
as the statement) of Lemma 4.

\medskip
 Step 1.
Let  $\mathbf{S}:=\overline{\{C_w: w\in W\}} \ \ \mbox{(the closure in $\mathbf{Aut}$)}$. 
Then $\mathbf{S}$ is a compact semigroup in $\mathbf{D}(\mu)$. Further, since $det C_w=e^{-qw}$, therefore it is an
invertible operator.  Thus
\[
\mathbf{H}:=\{C_w: w\in W\}\cup\{C_w^{-1}:w \in W\} \ \ \mbox{is  a commutative group in} \  \mathbf{Aut}.
\]
To this end we have check that for
 $w,u \in W$,  both $C_w\,C_u^{-1}$ and $C_w^{-1}C_u$ are in $\bf H$.   Let assume  that $w>u$ then  $C_wC_u^{-1}=C_{w-u}C_uC_u^{-1}=C_{w-u}\in \mathbf{H}$. Similarly, 
$C_w^{-1}C_u = (C_u C_{w-u})^{-1}C_u=C_{w-u}^{-1}C_u^{-1}C_u=C_{w-u}^{-1} \in \mathbf{H}$.
(These equations yield both closure and, with a trivial
extra argument, commutativity.)

\medskip
Step 2. Let  $ \bf G:= \mathbf{S}\cup \mathbf{S}^{-1}$. Then $\mathbf{G}\subset \mathbf{A}ut$ is a commutative compactly generated subgroup.
Moreover,  the mapping $ h: \mathbf{G}\to ( \Rset, +) $ given by $ h(A):=\log det A$ is a homomorphism of those two topological groups with the kernel $ \ker h= \mathbf{S}_0:=\mathbf{S}\cap \mathbf{A}(\mu)$.
Thus the quotient group $\mathbf{G}/ker h$ is isomorphic with $(\Rset, +)$.

\medskip
To see the above claim, first of all note that  since $\mathbf{S}_0$  is closed subsemigroup in the compact group $\mathbf{A}(\mu)$
therefore $\mathbf{S}_0$  is a compact group, by Theorem 1.1.12 in  Paalman - de Miranda (1964) (see Theorem 2
in the Appendix).

If $A\in \mathbf{S}_0$ then   $A\in \mathbf{A}(\mu)$  and by  Corollaries 2.3.2 and 2.4.2 from Jurek and Mason (1993), we have that  $|det A|=1$ . On the other hand, since $A\in \mathbf{S}$ we have that $0< detA\le1$, so $ det A=1$ and $h(A)=0$ and $\mathbf{S}_0\subset \ker h$. 

Conversely, if $det A=1$ and $A\in \mathbf{S}$  then $A\in\mathbf{D}(\mu)$ and by Jurek-Mason (1993), Proposition 2.3.5 we get that $A \in\mathbf{A}(\mu)$. Consequently, $A\in\mathbf{S}_0$. If $ det A=1$ and $A\in \mathbf{S}^{-1}$ then $A^{-1}\in\mathbf{S}$ and $ detA^{-1}=1$
so $A\in\mathbf{S}_0$, This completes the proof of the Step 2.

\medskip
Step 3.   There is an isomorphism $g:\mathbf{G}\to \Rset\oplus \mathbf{S}_0$ between the two topological groups.  

 This is so, because $\mathbf{G}$   is commutative and compactly generated group and the Pontriagin Theorem, from Montgomery and Zippin (1955), p.\ 187 (see Theorem 4
in the Appendix),  gives the needed isomorphism.

\medskip
Step 4. Taking the unit $\mathcal{I}$ in the group in $\mathbf{S}_0$ and putting for $t \geq 0$, 
\begin{multline}
T_t:=  g^{-1}(-t\oplus \mathcal{I}),  \ \ \ \mbox{if} \ \ \ \   g(\mathbf{S}) = (-\infty,0]\oplus \mathbf{S}_0   \\
T_t:=  g^{-1}(t\oplus \mathcal{I}),  \ \ \ \mbox{if} \ \ \ \   g(\mathbf{S}) = [0, \infty) \oplus \mathbf{S}_0 
\end{multline}
we obtain the one-parameter semigroup of matrices in  $\mathbf{D}(\mu)$.

>From the  equality
$g(\mathbf{G})=g(\mathbf{S})\cup (g(\mathbf{S}))^{-1}$,
 and the fact  $g(\mathbf{S})$ is closed subsemigroup we infer  that   either $g(\mathbf{S}) = (-\infty,0]\oplus \mathbf{S}_0 $ or $g(\mathbf{S}) = [0, \infty)\oplus \mathbf{S}_0$. 

\medskip
Step 5. For $t\ge 0$, $T_t=\exp(-t V)$ for some matrix $V$,  and $T_t\to 0$ as $t\to \infty$. 

\medskip
By Hille (1948), Theorem 8.4.2 
(or Hille and Phillips (1957), Theorem 9.4.2 --- see
Theorem 3 in the Appendix ---
with the idempotent there being the identity matrix here
in our context),  
we get the exponential form, that is, $T_t=\exp tQ, t\ge 0$, for some matrix $Q$. 

For $t>0$ we have that $T_t\notin\mathbf{S}_0$
and thus 
\begin{equation}
0<\det T_t<1 \  \  \mbox{for all} \ \ t>0
\end{equation}
From the definitions of operators $C_w, T_{n,r}$ and semigroup $\mathbf{S}$ it follows that the idempotents $J_r, 1\le r\le q$ commute with $\mathbf{S}$

Since $T_t\in \mathbf{D}(\mu), t \ge 0$, the set
$\{T_t, t \ge 0\}$ is conditionally compact. 
Hence by the Numakura Theorem, among the limits points (as $t\to\infty$) there is an idempotent, say $K$. 
Of course by (36) and a simple argument, 
$det K=0$; and by (24), $K=J_1K+...+J_qK$.
Also, $K$ is the limit of a sequence of $C_w$'s with
$w \to \infty$ (forced by (32) since $det K = 0$), and 
hence by (31),
\[
det _{J_r}K=0  \ \ \  (r=1,2,...,q)
\]
Since $K$ and $J_r$ commute and both are idempotents then so is $J_rK$. From above and (17),  $det_{J_r}J_rK=det_{J_r}K=0$, so $J_r\neq J_rK$. Moreover, we also have that $J_r(J_rK)=(J_rK)J_r=J_rK$. Thus  from the properties of $J_r$  
((24) and the {\it entire sentence\/} containing (25)) 
we must have $J_rK=0$ and consequently $K=0$. 
That is, the only limit point of $T_t$ as $t\to\infty$. 
As a consequence, Lemma 5 holds.

\medskip
\emph{Proof of Theorem 1.}
It follows from Lemma 5.

\medskip
\textbf {4. Appendix.}

For an ease of reference let us quote here the following algebraic facts.

\begin{thm}
 Each locally compact subsemigroup S of a compact group G is a compact subgroup.
\end{thm}
Cf.\ A.B.\ Paalman - De Miranda (1964), Theorem 1.1.12.

\medskip
\begin{thm}
If  $T:(0,\infty)\to \mathfrak{B} \  \mbox{(a real or complex Banach algebra)}$ \ \ \mbox{satisfies} 
\begin{equation*}
T (t+s)=T(t)T(s) \ \ \mbox{for all} \ \  0< t, s < \infty \  \ \  \mbox{and} \ \ \  
\lim_{t\to 0} T(t)=J \ \mbox{(an idempotent)},
\end{equation*}
then there exists an element $ A \in \mathfrak{B} $ such that
\begin{equation}
T(t)=J+\sum_{n=1}^{\infty}\frac{t^n}{n!}\, A^n \ \ \ \mbox{(absolutely convergent series)}.
\end{equation}
\end{thm}
Cf.\ E.\ Hille (1948), Theorem 8.4.2 or E.\ Hille and R.\ Phillips (1957), Theorem 9.4.2.

\medskip
\begin{thm} \emph{(Pontriagin Theorem)}  Suppose a topological group  $G^{\prime}$, generated  by a compact set, contains a compact subgroup $H^{\prime}$
such that  $G^{\prime}/H^{\prime}$ is isomorphic with an 
$r$-dimensional real vector group $V_r$. Then  $G^{\prime}$ has a vector subgroup $E_r$ such that 
$G^{\prime}= H^{\prime}\oplus E_r$
\end{thm}
Cf.\ Montgomery and Zippin (1955), p. 187.

\medskip
\medskip

\begin{center}
REFERENCES
\end{center}

\noindent [1]  A. Araujo and E. Gine (1980), \emph{The central limit theorem for real
and Banach valued random variables}, John Wiley $\&$ Sons, New York.

\noindent [2] P. Billingsley (1999), \emph{Convergence of probability measures},
Second edition, John Wiley $\&$ Sons, New York.

\noindent [3] R.C. Bradley (2007), \emph{Introduction to strong mixing conditions},  vol. 1 and  vol. 2, 
Kendrick Press, Heber City (Utah).

\noindent [4] R.C. Bradley and Z. J. Jurek (2014), The strong mixing and
the selfdecomposability properties, 
\emph {Statist. Probab. Letters} \textbf{84}, pp. 67-71.

\noindent [5] R.M. Dudley (2002), \emph{Real analysis and probability}, Second edition, 
Cambridge University Press, Cambridge.

\noindent [6] M.G. Hahn and M. Klass (1981), The multidimensional central limit theorem for arrays normed by affine transformations, \emph{Ann. Probab.} \textbf{9}, pp.\ 611-623.

\noindent [7] E. Hille (1948), \emph{Functional analysis and semi-groups,} American Mathematical Society, New York.

\noindent [8]  E. Hille and  R.S. Philips (1957), \emph{Functional analysis and semi-groups}, American Mathematical Society, Providence.

\noindent [9] Z.J. Jurek (1982),  An integral representation of operator-selfdecomposable random variables,
\emph{Bull. Acad. Polon. Sci.} \textbf{30}, pp. 385-393.

\noindent [10]  Z.J. Jurek  and J. David Mason (1993),  \emph{Operator-limit distributions in probability theory},  John Wiley $\&$ Sons, New York.

\noindent [11] D. Montgomery and L. Zippin (1955),  \emph{Topological transformation  groups},
Interscience Publishers,
New York and London.

\noindent [12]  K. Numakura (1952), On bicompact semigroups, 
\emph{ Math. J. Okayama University}, 
\textbf{1}, pp. 99-108.

\noindent [13]  A.B. Paalman-De Miranda (1964), \emph{Topological semigroups}, Mathematisch Centrum, Amsterdam.

\noindent [14]  K. R. Parthasarathy (1967), \emph{Probability measures on metric spaces}, 
Academic Press, New York and London.

\noindent [15] M. Rosenblatt (1956), A central limit theorem and a strong
mixing, \emph{Proc. Natl. Acad. Sci. USA} \textbf{42}, pp.
43-47.

\noindent [16]  M. Sharpe (1969), Operator-stable probability measures on vector groups, \emph{Trans. Amer. Math. Soc.} \textbf{136},  pp. 51-65.

\noindent [17]  K. Urbanik (1972),  L\'evy's probability measures on Euclidean spaces, 
\emph{Studia Math.} \textbf{44}, pp. 119-148.

\noindent [18]  K. Urbanik (1978), L\'evy's probability measures on Banach spaces, 
\emph{Studia Math.} \textbf{63}, pp. 283-308.

\end{document}